\documentclass[11pt,reqno]{amsart}
\usepackage[margin=1in]{geometry}
\usepackage{amsmath,amssymb,amsthm,mathrsfs}
\usepackage{microtype}
\usepackage{tikz}
\usetikzlibrary{arrows.meta}
\definecolor{pathblue}{HTML}{126677}
\definecolor{pathrust}{HTML}{A04A28}
\definecolor{pathink}{HTML}{28383F}
\definecolor{pathguide}{HTML}{879299}
\tikzset{
  time line/.style={draw=pathguide, line width=.45pt},
  guide/.style={draw=pathguide, densely dotted, line width=.5pt},
  left path/.style={draw=pathblue, line width=1.2pt, line cap=round, line join=round},
  right path/.style={draw=pathrust, line width=1.2pt, dash pattern=on 3pt off 2pt,
                     line cap=round, line join=round},
  common path/.style={draw=pathink, line width=1.9pt, line cap=round, line join=round},
  dimension/.style={{Stealth[length=3.2pt]}-{Stealth[length=3.2pt]},
                    draw=pathink, line width=.55pt},
  panel title/.style={font=\small\bfseries, anchor=west, inner sep=0pt},
  time label/.style={font=\scriptsize, anchor=east, inner sep=2pt},
  point label/.style={font=\small, inner sep=2pt}
}
\usepackage{hyperref}
\hypersetup{pdfauthor={Ruixuan Zhang},pdftitle={Coalescence and the Two-Point Covariance in the Directed Landscape},
  colorlinks=true,linkcolor=blue,citecolor=blue,urlcolor=blue}

\newtheorem{theorem}{Theorem}[section]
\newtheorem{lemma}[theorem]{Lemma}
\newtheorem{proposition}[theorem]{Proposition}
\newtheorem{corollary}[theorem]{Corollary}
\theoremstyle{remark}
\newtheorem{remark}[theorem]{Remark}
\numberwithin{equation}{section}
\newcommand{\R}{\mathbb{R}}
\newcommand{\E}{\mathbb{E}}
\renewcommand{\P}{\mathbb{P}}
\newcommand{\Var}{\operatorname{Var}}
\newcommand{\Cov}{\operatorname{Cov}}
\newcommand{\sgn}{\operatorname{sgn}}
\newcommand{\ind}{\mathbf{1}}

\begin{document}

\title[Coalescence and Covariance in the Directed Landscape]{Coalescence and the Two-Point Covariance in the Directed Landscape}

\author{Ruixuan Zhang}
\address{Department of Mathematics, University of Utah, 155 S 1400 E, JWB 317,
  Salt Lake City, UT 84112, USA}
\email{ray.zhang@math.utah.edu}

\begin{abstract}
In the directed landscape with Brownian terminal data, we express the spatial derivative of the two-point height covariance in terms of the probability that two geodesics coalesce before the terminal time and the exit law of a single geodesic. This formula arises as the Kardar--Parisi--Zhang scaling limit of a zero-temperature Brownian last-passage identity from our earlier work on the O'Connell--Yor polymer. To pass the coalescence probabilities to the limit, we show that two paths sharing an optimal exit almost surely lose total weight if they must stay disjoint before reaching it. At stationarity, we determine when two zero-direction semi-infinite geodesics first coalesce. We recover this distribution by differentiating the height-increment variance first in the starting separation, then in the boundary strength. When the terminal reward is zero everywhere, the starting points with the same leftmost optimal exit as the origin form an interval. The probability that this interval extends past given distances on both sides equals the rate at which the height covariance decreases at their sum.
\end{abstract}

\maketitle

\section{Introduction}

The directed landscape is a random directed metric on space-time that arises in the Kardar--Parisi--Zhang (KPZ) universality class. For $x,y\in\R$ and $s<t$, its value $\mathcal L(x,s;y,t)$ is the optimal passage weight from $(x,s)$ to $(y,t)$. At any intermediate time, this weight is the maximum, over all crossing locations, of the sum of the two passage weights \cite[Definition~10.1]{DOV22}. Dauvergne, Ortmann and Vir\'ag constructed the landscape as the scaling limit of Brownian last-passage percolation (BLPP). On the scale of $N$ layers, transverse displacements are of order $N^{2/3}$ and centered passage times are of order $N^{1/3}$ \cite[Theorem~1.5]{DOV22}. BLPP itself also arises as a diffusive limit of lattice last-passage times with a fixed number of rows \cite[Theorems~3.1--3.2]{GW91}.

The O'Connell--Yor polymer is a semi-discrete directed random polymer in an environment of independent Brownian motions, one per layer. The layers are discrete, but a path can jump between them at continuous positions. Its energy $E$ is the sum of the Brownian increments it collects. For fixed endpoints, the partition function $Z$ integrates $e^{\beta E}$ over the admissible jump positions \cite[Section~5]{OY01}. At fixed inverse temperature $\beta>0$, the free energy $\beta^{-1}\log Z$ converges in law to the directed landscape after centering and KPZ rescaling. The convergence is joint for any finite collection of time pairs, in the topology of locally uniform convergence in the spatial endpoints \cite[Theorem~7.4]{DZ24polymer}. Fix the number of layers and the horizontal endpoints $X<P$. As $\beta\to\infty$, the free energy converges to the largest weight among these paths, the BLPP passage time.

In earlier work with Rassoul-Agha, Shen and Zheng, we used Gaussian integration by parts to obtain an identity for the free energy $F_{N,\beta}$ of the $N$-layer polymer with Brownian terminal reward of strength $\alpha$. For unscaled starts $a<b$,
\[
 -\partial_b\Cov(F_{N,\beta}(a),F_{N,\beta}(b))
 =\tfrac12(1-\alpha^2)\P_\beta(M)+\alpha^2\P_\beta(Q_b\le0).
\]
Here $M$ denotes a meeting before layer $N$ and $Q_b$ is the exit of the right path. The probability $\P_\beta$ averages over the environment and two paths sampled independently given that environment \cite[Theorem~1.1, (1.4)--(1.5)]{RSZZ26}. At $\beta=\infty$ the same identity holds for BLPP passage values and their maximizing paths \cite[Corollary~1.12]{RSZZ26}. We take the $N\to\infty$ KPZ scaling limit of this zero-temperature formula.

At time $T>0$ a path may exit anywhere on the line and receives the terminal reward $f(y)=\sqrt2\alpha B(y)$, where $B$ is a two-sided standard Brownian motion pinned at zero and independent of $\mathcal L$. The optimal height from $(x,0)$ is
\begin{equation}\label{eq:height}
 H(x)=\max_{y\in\R}\{\mathcal L(x,0;y,T)+f(y)\}.
\end{equation}
All starts share the same $\mathcal L$ and $B$. For each fixed $x$, the maximizing exit $Y(x)$ and its geodesic $\pi_x$ are almost surely unique (Lemmas~\ref{lem:heights_exit} and~\ref{lem:free_geodesics}). For fixed $a<b$, the event of an interior meeting is
\[
 M(a,b)=\{\pi_a(r)=\pi_b(r)\text{ for some }0<r<T\}.
\]
\begin{theorem}[Covariance and coalescence]\label{thm:main}
The function $(a,b)\mapsto\Cov(H(a),H(b))$ is continuously differentiable on $a<b$, with
\begin{equation}\label{eq:cov_deriv}
 \begin{aligned}
 -\partial_b\Cov(H(a),H(b))
  &=(1-\alpha^2)\P(M(a,b))+2\alpha^2\P(Y(b)\le 0),\\
 \partial_a\Cov(H(a),H(b))
  &=(1-\alpha^2)\P(M(a,b))+2\alpha^2\P(Y(a)>0).
 \end{aligned}
\end{equation}
Moreover, $M(a,b)=\{Y(a)=Y(b)\}$ almost surely.
\end{theorem}

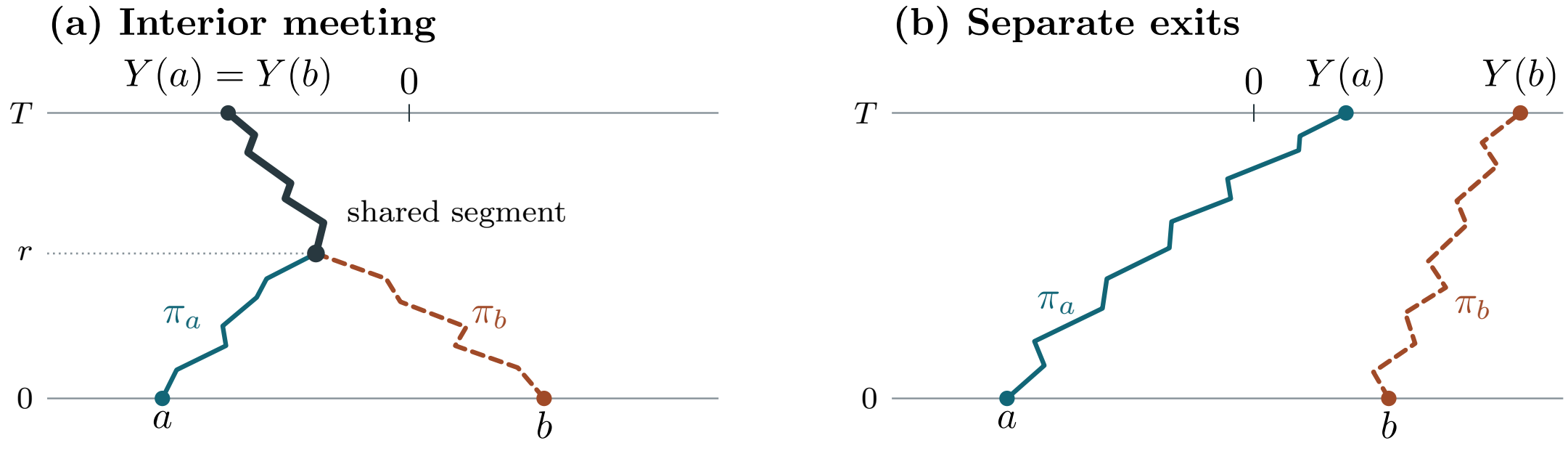
\begin{figure}[htbp]
\centering
\begin{tikzpicture}[x=1cm,y=1cm,font=\small]
  \begin{scope}
    \node[panel title] at (0,3.40) {(a) Interior meeting};
    \draw[time line] (0,0)--(6.12,0) (0,2.60)--(6.12,2.60);
    \node[time label] at (-.05,0) {$0$};
    \node[time label] at (-.05,2.60) {$T$};
    \draw[guide] (0,1.32)--(2.45,1.32);
    \node[time label] at (-.05,1.32) {$r$};
    \draw[left path] (1.05,0)--(1.18,.26)--(1.63,.48)--(1.60,.66)
      --(1.91,.92)--(2.00,1.09)--(2.45,1.32);
    \draw[right path] (4.53,0)--(4.29,.28)--(3.72,.48)--(3.82,.65)
      --(3.22,.88)--(3.09,1.09)--(2.45,1.32);
    \draw[common path] (2.45,1.32)--(2.52,1.60)--(2.17,1.82)
      --(2.22,1.96)--(1.83,2.24)--(1.89,2.40)--(1.65,2.60);
    \fill[pathblue] (1.05,0) circle (2pt);
    \fill[pathrust] (4.53,0) circle (2pt);
    \fill[pathink] (2.45,1.32) circle (2.3pt) (1.65,2.60) circle (2pt);
    \draw[pathink] (3.30,2.52)--(3.30,2.68);
    \node[point label,above] at (3.30,2.70) {$0$};
    \node[point label,above] at (1.65,2.68) {$Y(a)=Y(b)$};
    \node[point label,below] at (1.05,-.04) {$a$};
    \node[point label,below] at (4.53,-.04) {$b$};
    \node[point label,left,text=pathblue] at (1.50,.74) {$\pi_a$};
    \node[point label,right,text=pathrust] at (3.79,.74) {$\pi_b$};
    \node[font=\scriptsize,anchor=west] at (2.60,1.68) {shared segment};
  \end{scope}
  \begin{scope}[xshift=7.70cm]
    \node[panel title] at (0,3.40) {(b) Separate exits};
    \draw[time line] (0,0)--(6.12,0) (0,2.60)--(6.12,2.60);
    \node[time label] at (-.05,0) {$0$};
    \node[time label] at (-.05,2.60) {$T$};
    \draw[left path] (1.05,0)--(1.39,.30)--(1.30,.52)--(1.92,.82)
      --(1.96,1.09)--(2.53,1.37)--(2.55,1.61)--(3.09,1.82)
      --(3.06,2.00)--(3.71,2.26)--(3.72,2.39)--(4.14,2.60);
    \draw[right path] (4.53,0)--(4.39,.24)--(4.77,.50)--(4.68,.78)
      --(5.05,1.01)--(4.89,1.25)--(5.24,1.58)--(5.15,1.80)
      --(5.52,2.12)--(5.38,2.33)--(5.73,2.60);
    \fill[pathblue] (1.05,0) circle (2pt) (4.14,2.60) circle (2pt);
    \fill[pathrust] (4.53,0) circle (2pt) (5.73,2.60) circle (2pt);
    \draw[pathink] (3.30,2.52)--(3.30,2.68);
    \node[point label,above] at (3.30,2.70) {$0$};
    \node[point label,above] at (4.14,2.68) {$Y(a)$};
    \node[point label,above] at (5.73,2.68) {$Y(b)$};
    \node[point label,below] at (1.05,-.04) {$a$};
    \node[point label,below] at (4.53,-.04) {$b$};
    \node[point label,left,text=pathblue] at (1.77,.87) {$\pi_a$};
    \node[point label,right,text=pathrust] at (5.05,.83) {$\pi_b$};
  \end{scope}
\end{tikzpicture}
\caption{Two configurations in the covariance identity \eqref{eq:cov_deriv}, with time increasing upward. The tick at $0$ on the terminal line marks the Brownian pin. In (a), the paths meet at time $r<T$ and share a negative exit. In (b), they stay disjoint and both exits are positive. The meeting probability carries weight $1-\alpha^2$ in both derivatives. The exit probability carries weight $2\alpha^2$ and uses the exit of the start that moves. For $-\partial_b$ the event is $Y(b)\le0$, as in (a). For $\partial_a$ it is $Y(a)>0$, as in (b).}
\label{fig:covariance-events}
\end{figure}

Figure~\ref{fig:covariance-events} illustrates the events in \eqref{eq:cov_deriv}. For flat data ($\alpha=0$), every exit receives zero terminal reward, and the rate of covariance loss as $b$ increases is the meeting probability. At $|\alpha|=1$ the rate is $2\P(Y(b)\le0)$, and $H(\cdot)-H(0)$ is a two-sided Brownian motion of variance rate $2$ \cite[Lemma~10.2(3)]{DOV22}, \cite[Theorem~4.5(ii)]{MQR21}. Integration in $b$ gives the covariance itself.

\begin{corollary}[Integrated covariance formula]\label{cor:integral}
For $a<b$,
\begin{equation}\label{eq:cov_integral}
 \Cov(H(a),H(b))=\Var(H(a))
 -(1-\alpha^2)\int_a^b\P(M(a,u))\,du
 -2\alpha^2\int_a^b\P(Y(u)\le 0)\,du .
\end{equation}
\end{corollary}

In Section~\ref{sec:finite} a linear tilt and edge estimates for Hermite and Wishart matrices give uniform fourth moments for exits and heights (Lemma~\ref{lem:heights_exit}). These bounds localize the maximizers on the unbounded terminal line and give uniform integrability for the covariance. Section~\ref{sec:coalescence} treats meetings. Disjoint paths may converge uniformly to paths that touch at the terminal time or share a whole segment. At a common maximizing exit, the two-line representation and local Brownian absolute continuity \cite[Theorems~1.2 and~1.4]{Dau23} give a strictly positive loss in total weight when the paths are required to stay disjoint before $T$. Joint convergence of the one- and two-path weights \cite[Theorems~1.5--1.6]{DZ24} then gives convergence of the meeting probabilities. We present the consequences of the identity first, in Section~\ref{sec:applications}.

\section{Discussion and applications}\label{sec:applications}

The choice of terminal reward determines which geodesic event can be read from the covariance. For flat data, the Airy$_1$ covariance determines the common-exit probability and the joint law of the two ends of the exit cluster containing the origin (Corollaries~\ref{cor:app-flat} and~\ref{cor:app-cluster}). For nearby starts, counting distinct exits gives the first-order probability of separation (Proposition~\ref{prop:app-short}). At stationary strength the meeting term drops out of the spatial derivative. Its coefficient changes with boundary strength, however, and this change recovers the distribution of the first coalescence time (Corollary~\ref{cor:app-stationary}).

Write $q=\alpha^2$ and use $\sqrt{2q}\,B$ for the terminal reward. We write $H_{q,T}$ and $Y_{q,T}$ when the boundary strength or time varies. For $r>0$, let $p_{q,T}(r)=\P(M(0,r))$ and set $p_{q,T}(0)=1$. Upon translating both starts by $a$, recentering $B$ at $a$ subtracts the same reward from every competitor. Spatial translation and $1{:}2{:}3$ scaling of the landscape \cite[Lemma~10.2]{DOV22}, together with Brownian scaling, give
\begin{equation}\label{eq:app-scaling}
 \P(M(a,a+r))=p_{q,T}(r)=p_{q,1}(rT^{-2/3}),\qquad
 Y_{q,T}(a)-a\ \overset d=\ T^{2/3}Y_{q,1}(0).
\end{equation}
For fixed parameters we continue to write $Y(x)$ for the exit. The proof for nearby starts uses the exit moment bound of Section~\ref{sec:finite} and uniqueness of the free geodesic (Lemma~\ref{lem:free_geodesics}). The stationary response formula also uses the strict gap at a common exit (Lemma~\ref{lem:strict_loss}).

\subsection{Flat terminal data and the \texorpdfstring{Airy$_1$}{Airy 1} covariance}
Let $C_T(r)=\Cov(H_{0,T}(0),H_{0,T}(r))$ and let $g_1(u)=\Cov(\mathcal A_1(0),\mathcal A_1(u))$ in the normalization of Matetski, Quastel and Remenik \cite[(4.15)]{MQR21}. Two starts far apart can select the same exit, even though the typical displacement of either exit is of order $T^{2/3}$. The flat covariance determines the probability of this event at every separation.

\begin{corollary}[Coalescence at a flat boundary]\label{cor:app-flat}
For $T,r>0$,
\begin{equation}\label{eq:app-airy}
 p_{0,T}(r)=-g_1'\!\left(\frac{r}{(2T)^{2/3}}\right),\qquad
 \lim_{r/T^{2/3}\to\infty}\frac{T^2}{r^3}\log p_{0,T}(r)=-\frac13 .
\end{equation}
\end{corollary}

\begin{proof}
Landscape flip symmetry \cite[Lemma~10.2(3)]{DOV22} and the flat height law \cite[(4.15)]{MQR21} give $C_T(r)=(2T)^{2/3}g_1(r/(2T)^{2/3})$. Theorem~\ref{thm:main} gives $-C_T'(r)=p_{0,T}(r)$, proving the first identity in \eqref{eq:app-airy}.

Basu, Busani and Ferrari \cite[Theorem~1.1]{BBF23} prove $\log g_1(u)\sim-4u^3/3$, so $\log C_1(r)\sim-r^3/3$. For $0<r<s$, Lemma~\ref{lem:free_geodesics} gives $Y(0)\le Y(r)\le Y(s)$. A common exit for $0$ and $s$ is therefore also a common exit for $0$ and $r$, so $p_{0,1}$ is nonincreasing. For $0<\varepsilon<1$, integrating $-C_1'=p_{0,1}$ gives
\[
 \frac{C_1(r)-C_1((1+\varepsilon)r)}{\varepsilon r}
 \le p_{0,1}(r)\le
 \frac{C_1((1-\varepsilon)r)-C_1(r)}{\varepsilon r}.
\]
For fixed $0<a<b$, the same covariance asymptotic gives
\[
 \log\frac{C_1(br)}{C_1(ar)}
 =-\frac{b^3-a^3}{3}\,r^3+o(r^3).
\]
Thus $C_1((1+\varepsilon)r)/C_1(r)$ and $C_1(r)/C_1((1-\varepsilon)r)$ both tend to zero. Taking logarithms of the two bounds and dividing by $r^3$ gives limits $-1/3$ and $-(1-\varepsilon)^3/3$, respectively. Let $\varepsilon\downarrow0$ and apply \eqref{eq:app-scaling}.
\end{proof}

\subsection{Common exits and the length of an exit cluster}
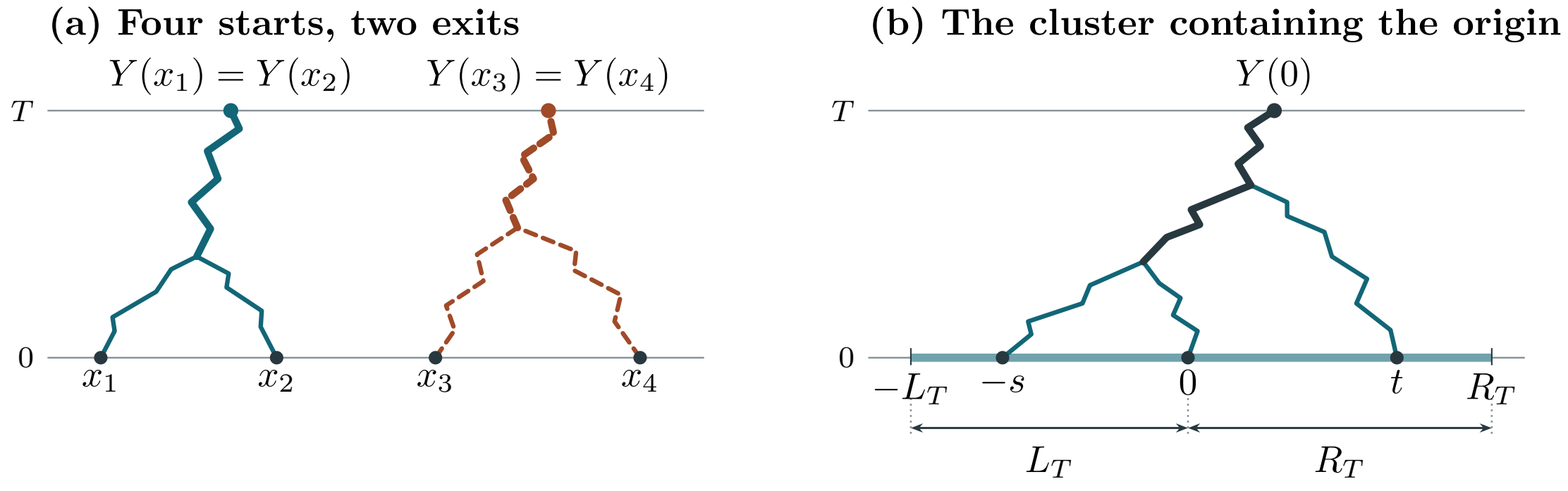
\begin{figure}[htbp]
\centering
\begin{tikzpicture}[x=1cm,y=1cm,font=\small]
  \begin{scope}
    \node[panel title] at (0,3.12) {(a) Four starts, two exits};
    \draw[time line] (0,0)--(6.16,0) (0,2.32)--(6.16,2.32);
    \node[time label] at (-.05,0) {$0$};
    \node[time label] at (-.05,2.32) {$T$};
    \draw[left path] (.50,0)--(.63,.25)--(.61,.38)--(1.02,.62)--(1.16,.83)
      --(1.40,.95);
    \draw[left path] (2.15,0)--(2.00,.29)--(2.01,.44)--(1.68,.66)
      --(1.70,.79)--(1.40,.95);
    \draw[left path,line width=1.9pt] (1.40,.95)--(1.53,1.21)--(1.35,1.46)
      --(1.60,1.68)--(1.50,1.94)--(1.80,2.15)--(1.72,2.32);
    \draw[right path] (3.64,0)--(3.82,.26)--(3.74,.49)--(4.09,.72)
      --(4.03,.97)--(4.41,1.22);
    \draw[right path] (5.56,0)--(5.32,.34)--(5.38,.59)--(4.94,.82)
      --(4.96,1.00)--(4.41,1.22);
    \draw[right path,line width=1.9pt] (4.41,1.22)--(4.30,1.48)--(4.56,1.68)
      --(4.45,1.89)--(4.75,2.10)--(4.70,2.32);
    \foreach \x/\name in {.50/1,2.15/2,3.64/3,5.56/4}{
      \fill[pathink] (\x,0) circle (1.8pt);
      \node[point label,below] at (\x,-.04) {$x_{\name}$};}
    \fill[pathblue] (1.72,2.32) circle (2pt);
    \fill[pathrust] (4.70,2.32) circle (2pt);
    \node[point label,above] at (1.72,2.39) {$Y(x_1)=Y(x_2)$};
    \node[point label,above] at (4.70,2.39) {$Y(x_3)=Y(x_4)$};
  \end{scope}
  \begin{scope}[xshift=7.70cm]
    \node[panel title] at (0,3.12) {(b) The cluster containing the origin};
    \draw[time line] (0,0)--(6.16,0) (0,2.32)--(6.16,2.32);
    \node[time label] at (-.05,0) {$0$};
    \node[time label] at (-.05,2.32) {$T$};
    \draw[draw=pathblue!60,line width=2.2pt] (.40,0)--(5.85,0);
    \foreach \x in {.4,3,5.85}{\draw[guide] (\x,-.38)--(\x,-.73);}
    \foreach \x in {.40,5.85}{\draw[pathink] (\x,-.07)--(\x,.09);}
    \draw[left path] (1.26,0)--(1.53,.22)--(1.50,.34)--(2.01,.51)
      --(2.08,.68)--(2.58,.90);
    \draw[left path] (3,0)--(3.09,.25)--(2.86,.40)--(2.93,.56)
      --(2.73,.70)--(2.58,.90);
    \draw[common path] (2.58,.90)--(2.80,1.13)--(3.11,1.25)--(3.03,1.39)
      --(3.59,1.62);
    \draw[left path] (4.96,0)--(4.90,.26)--(4.58,.51)--(4.68,.74)
      --(4.35,.95)--(4.28,1.18)--(3.93,1.33)--(3.93,1.46)--(3.59,1.62);
    \draw[common path] (3.59,1.62)--(3.47,1.82)--(3.68,1.99)--(3.56,2.16)
      --(3.81,2.32);
    \fill[pathink] (1.26,0) circle (1.8pt) (3,0) circle (1.8pt)
      (4.96,0) circle (1.8pt) (3.81,2.32) circle (2pt);
    \node[point label,above] at (3.81,2.39) {$Y(0)$};
    \node[point label,below] at (.4,-.07) {$-L_T$};
    \node[point label,below] at (1.26,-.04) {$-s$};
    \node[point label,below] at (3,-.04) {$0$};
    \node[point label,below] at (4.96,-.04) {$t$};
    \node[point label,below] at (5.85,-.07) {$R_T$};
    \draw[dimension] (.4,-.66)--(3,-.66)
      node[midway,below=1pt] {$L_T$};
    \draw[dimension] (3,-.66)--(5.85,-.66)
      node[midway,below=1pt] {$R_T$};
  \end{scope}
\end{tikzpicture}
\caption{Common exits and cluster length. In (a), only the step from $x_2$ to $x_3$ changes the exit, so $K=2$. In (b), the highlighted interval lies on the starting line and has length $L_T+R_T$. The test starts $-s$ and $t$ lie strictly inside its ends, so $Y(-s)=Y(0)=Y(t)$. Their separation is $s+t$, the argument of the coalescence probability in \eqref{eq:app-cluster}.}
\label{fig:exit-cluster}
\end{figure}

For ordered starts, the outermost two paths determine whether all paths share a terminal segment. In the flat case this also gives the law of the interval of starts that share an exit with the origin. Figure~\ref{fig:exit-cluster} shows the exit count for a finite family and the two distances that determine the length of this interval.

\begin{corollary}[Common exits and cluster lengths]\label{cor:app-cluster}
For fixed $x_1<\cdots<x_k$, $k\ge2$, let $K=\#\{Y(x_1),\ldots,Y(x_k)\}$. The $k$ paths share a terminal segment of positive duration precisely on $\{K=1\}$, almost surely. Their exit count satisfies
\begin{align}
 \P(K=1)&=p_{q,T}(x_k-x_1),\label{eq:app-all}\\
 \E K&=1+\sum_{i=2}^k\bigl[1-p_{q,T}(x_i-x_{i-1})\bigr].\label{eq:app-count}
\end{align}
At $q=0$, choose the leftmost maximizing exit $Y(x)$ at each real start. The cluster $I_T=\{x:Y(x)=Y(0)\}$ is an interval with ends $-L_T=\inf I_T$ and $R_T=\sup I_T$. For $s,t\ge0$ with $s+t>0$,
\begin{equation}\label{eq:app-cluster}
 \P(L_T>s,\ R_T>t)=p_{0,T}(s+t)=-C_T'(s+t).
\end{equation}
In particular $C_T''(dr)=\P(R_T\in dr)$ as measures on $(0,\infty)$, and
$\E(L_T+R_T)=2\Var(H_{0,T}(0))$.
\end{corollary}

\begin{proof}
The exit ordering in Lemma~\ref{lem:free_geodesics} gives
\[
 \{Y(x_1)=Y(x_k)\}=\{K=1\},\qquad
 K=1+\sum_{i=2}^k\ind_{\{Y(x_{i-1})<Y(x_i)\}}.
\]
On $\{K=1\}$ each of the finitely many paths has merged with the first strictly before $T$ by Theorem~\ref{thm:main}. After the last such meeting they all agree. This proves \eqref{eq:app-all} and \eqref{eq:app-count}.

At $q=0$ the parabolic bound \cite[Corollary~10.7]{DOV22} confines all maximizing exits from $[-n,n]$ to one random compact interval, as in the proof of Lemma~\ref{lem:heights_exit}. Taking $n\in\mathbb N$ gives a leftmost maximizing exit simultaneously for every real start. If $x<x'$ and $u=Y(x)>v=Y(x')$, take geodesics from $x$ to $u$ and from $x'$ to $v$ \cite[Theorem~1.7]{DZ24}. Their spatial difference is continuous and has opposite signs at $0$ and $T$. They meet at a time in $(0,T)$, where exchanging tails gives
\[
 H(x)+H(x')\le\mathcal L(x,0;v,T)+\mathcal L(x',0;u,T)\le H(x)+H(x').
\]
Thus $v<u$ also maximizes at $x$, a contradiction. Hence $Y$ is nondecreasing. Its level set $I_T$ is an interval containing $0$ and has the same infimum and supremum as $I_T\cap\mathbb Q$. For $\varepsilon>0$, exit order and translation give
\[
 p_{0,T}(s+t+2\varepsilon)\le \P(L_T>s,\ R_T>t)\le p_{0,T}(s+t).
\]
Theorem~\ref{thm:main} gives continuity at positive separation, so letting $\varepsilon\downarrow0$ proves \eqref{eq:app-cluster}. The one-sided bounds give $\P(R_T>r)=\P(L_T>r)=-C_T'(r)$ for $r>0$. Taking the distributional derivative gives the measure identity on $(0,\infty)$. At the origin, Proposition~\ref{prop:app-short} gives $p_{0,T}(r)\to1$ as $r\downarrow0$. Taking $s=t=r$ in \eqref{eq:app-cluster} yields
\[
 \P(L_T>0,\ R_T>0)=\lim_{r\downarrow0}p_{0,T}(2r)=1.
\]
Thus neither endpoint distance has an atom at zero, and the origin lies in the interior of $I_T$ almost surely. At $q=0$, \eqref{eq:cov_integral} gives $C_T(0)-C_T(r)=\int_0^r p_{0,T}(u)\,du$, so $C_T(0+)=C_T(0)=\Var(H_{0,T}(0))$. Since $C_T(r)\to0$ by \cite[Theorem~1.1]{BBF23},
\[
 \E L_T=\E R_T
 =\int_0^\infty\P(R_T>r)\,dr
 =\int_0^\infty -C_T'(r)\,dr=C_T(0).
\]
\end{proof}

\subsection{Nearby starts}
The number of distinct exits on a fine grid determines the first-order probability of separation. The coefficient below is the expected number of exit changes across a unit interval at time one.

\begin{proposition}[Separation of nearby starts]\label{prop:app-short}
For every fixed $q\ge0$ there is a constant $\lambda_q\in(0,\infty)$ such that
\begin{equation}\label{eq:app-short}
 1-p_{q,T}(r)\sim\lambda_q\,rT^{-2/3}
 \qquad\text{as }rT^{-2/3}\downarrow0.
\end{equation}
At $T=1$, if $K_m$ counts the distinct exits from the grid $2^{-m}\mathbb Z\cap[0,1]$, then $\lambda_q=\E[\lim_{m\to\infty}K_m-1]$.
\end{proposition}

\begin{proof}
Take $T=1$ and put $F(r)=1-p_{q,1}(r)$. For $r\ge0$ and $h>0$, translation in \eqref{eq:app-scaling} and the ordering $Y(-r)\le Y(0)\le Y(h)$ from Lemma~\ref{lem:free_geodesics} give
\[
 \begin{aligned}
 F(r+h)-F(r)
 &=\P\{Y(-r)<Y(h)\}-\P\{Y(-r)<Y(0)\}\\
 &=\P\{Y(-r)=Y(0)<Y(h)\}.
 \end{aligned}
\]
As $r$ increases, the event on the right decreases by the same ordering. Thus $F$ is concave, and $F(0)=0$ implies that $F(r)/r$ is nonincreasing on $(0,\infty)$. By \eqref{eq:app-count},
\[
 2^m F(2^{-m})=\E K_m-1,\qquad K_m\uparrow K_\infty.
\]
It remains to bound $\E K_\infty$. The union of the dyadic grids is countable, so Lemma~\ref{lem:free_geodesics} holds simultaneously for all pairs of grid starts. In particular every grid exit lies in $[Y(0),Y(1)]$. For $a=0,1$, Lemma~\ref{lem:heights_exit} gives $Y_N(a)\Rightarrow Y(a)$, so Fatou's lemma in a Skorokhod coupling and \eqref{eq:heights_exit_moment} give
\[
 \E|Y(a)|^4\le\liminf_{N\to\infty}\E|Y_N(a)|^4<\infty .
\]
A free optimizer is the unique geodesic between its realized endpoints, since any alternative would also maximize the free problem. Paths with different exits occupy different positions at time $1/2$, by \eqref{eq:tails}. Let $J_\ell$ count the distinct midpoints of unique geodesics with both spatial endpoints in $[-\ell,\ell]$ at times $0,1$. Ganguly and Zhang \cite[Lemma~3.12]{GZ} give $\P(J_\ell>M)\le Ce^{-cM^{1/384}}$ for $M>\ell^3$. Take $\ell=M^{3/10}$ for $M>1$. On $\{\max(|Y(0)|,|Y(1)|)\le\ell\}$ every grid start and exit lies in $[-\ell,\ell]$. Distinct exits therefore contribute distinct midpoints to $J_\ell$, so $K_\infty\le J_\ell$ on this event. Since $M>\ell^3$,
\[
 \P(K_\infty>M)
 \le \P\bigl(\max\{|Y(0)|,|Y(1)|\}>\ell\bigr)+\P(J_\ell>M)
 \le C_q M^{-6/5}+Ce^{-cM^{1/384}} .
\]
The tail bound is summable, so $\E K_\infty<\infty$. By concavity, the limit of $F(r)/r$ as $r\downarrow0$ agrees with its dyadic limit, and monotone convergence gives
\[
 \lim_{r\downarrow0}\frac{F(r)}{r}
 =\lim_{m\to\infty}(\E K_m-1)=\E K_\infty-1.
\]
If this limit were zero, $K_\infty=1$ would force $Y(0)=Y(1)$ almost surely. The exits are integrable, however, and their translation law in \eqref{eq:app-scaling} gives $\E Y(1)-\E Y(0)=1$. The time scaling in the same formula gives the result for every $T>0$.
\end{proof}

\subsection{The stationary response and the first coalescence time}
At $q=1$ the coefficient of the meeting probability in \eqref{eq:cov_deriv} vanishes. Varying the boundary strength recovers that probability. For $r>0$ let
\[
 V_{q,T}(r)=\Var\bigl(H_{q,T}(r)-H_{q,T}(0)\bigr).
\]
Let $\tau_r$ be the first coalescence time of the zero-direction semi-infinite geodesics from $(0,0)$ and $(r,0)$ in the same landscape.

\begin{corollary}[Stationary response and coalescence time]\label{cor:app-stationary}
For $q\ge0$ and $T,r>0$,
\begin{equation}\label{eq:app-variogram}
 \partial_r V_{q,T}(r)=2q+2(1-q)p_{q,T}(r).
\end{equation}
At stationary strength,
\begin{equation}\label{eq:app-time}
 \P(\tau_r<T)=\P(\tau_r\le T)=p_{1,T}(r)
 =1-\frac12\left.\partial_q\partial_r V_{q,T}(r)\right|_{q=1}.
\end{equation}
In particular, for fixed $r>0$, $\P(\tau_r>T)\sim\lambda_1 r T^{-2/3}$ as $T\to\infty$.
\end{corollary}

\begin{proof}
\textbf{Step 1. The spatial derivative.}
Write $v(x)=\Var(H_{q,T}(x))$. Pimentel \cite[Theorem~2, (1.6)]{Pim22}, with boundary variance $2q$, time reversal and \eqref{eq:app-scaling}, gives
$v'(r)=2q[2\P(Y(r)>0)-1]$ for $q>0$. At $q=0$ the flat height field is stationary, so $v'=0$. Expanding the variance and applying \eqref{eq:cov_deriv} gives
\[
 \partial_r V_{q,T}(r)
 =v'(r)+2(1-q)p_{q,T}(r)+4q\P(Y(r)\le0)
 =2q+2(1-q)p_{q,T}(r).
\]

\noindent\textbf{Step 2. The derivative in boundary strength.}
Take $q_n\to1$ and use the same $\mathcal L,B$ for every $q_n$ and for $q=1$. For $x\in\{0,r\}$ the parabolic bound in the proof of Lemma~\ref{lem:heights_exit} confines all these exits to one random interval $[-R,R]$. Almost surely the objectives $\mathcal L(x,0;y,T)+\sqrt{2q_n}\,B(y)$ converge uniformly there, since
\[
 |\sqrt{2q_n}-\sqrt2|\sup_{|y|\le R}|B(y)|\longrightarrow0 .
\]
Every subsequential limit of the maximizing exits is therefore a maximizer at $q=1$. By uniqueness for each of the two starts,
\[
 (u_n,v_n):=(Y_{q_n,T}(0),Y_{q_n,T}(r))
       \longrightarrow (u,v):=(Y_{1,T}(0),Y_{1,T}(r))\quad\text{a.s.}
\]
If $u<v$, eventually $u_n<v_n$. If $u=v$, the loss in total weight on requiring disjoint paths is strictly positive at $(u,v)$ by Lemma~\ref{lem:strict_loss}. The one- and two-path fields are continuous in the endpoints, so the loss at $(u_n,v_n)$ is eventually positive. Two disjoint free optimizers would attain zero loss. Thus in both cases the meeting indicators converge almost surely. Taking expectations gives $p_{q_n,T}(r)\to p_{1,T}(r)$. Dividing
\[
 \partial_r V_{q,T}(r)-2=2(q-1)\bigl[1-p_{q,T}(r)\bigr]
\]
by $q-1$ and letting $q\to1$ gives the derivative in \eqref{eq:app-time}.

\noindent\textbf{Step 3. The first coalescence time.}
At time $T$ the zero-direction Busemann profile, recentered at zero, is a two-sided Brownian motion of variance rate $2$, independent of the landscape on $[0,T]$ \cite[Theorem~3.13]{GSS26}. Use this profile as the common terminal reward. For $x=0,r$, the semi-infinite geodesic from $(x,0)$ passes through $(Y_{1,T}(x),T)$ by \cite[(23)]{GSS26}. By metric composition its restriction to $[0,T]$ is a free optimizer, unique by Lemma~\ref{lem:free_geodesics}. The two deterministic starts have unique zero-direction geodesics, and these agree from their first meeting onward \cite[Corollary~3.20(2)--(3)]{RV25}. Thus, almost surely,
\[
 \{\tau_r\le T\}=\{Y_{1,T}(0)=Y_{1,T}(r)\}=M(0,r)=\{\tau_r<T\},
\]
where the middle equality is Theorem~\ref{thm:main}. Taking probabilities gives \eqref{eq:app-time}, and Proposition~\ref{prop:app-short} gives the tail.
\end{proof}

\begin{remark}[Covariance of disjoint increments]\label{rem:increments}
For $a<b<c<d$ and $c<x<d$, take the difference of the derivatives with respect to the right starting point $x$ in \eqref{eq:cov_deriv}, evaluated at $(b,x)$ and $(a,x)$. The single-exit terms cancel. Integration over $c<x<d$ and the ordering $Y(a)\le Y(b)\le Y(x)$ from Lemma~\ref{lem:free_geodesics} give
\begin{equation}\label{eq:app-increments}
 \Cov\bigl(H(b)-H(a),H(d)-H(c)\bigr)
 =(q-1)\int_c^d \P\{Y(a)<Y(b)=Y(x)\}\,dx .
\end{equation}
Thus the increment covariance is nonpositive for $q<1$ and nonnegative for $q>1$. Its magnitude is determined by the probability that the two rightmost starts share an exit and the leftmost start has a separate exit.
\end{remark}

\subsection{Connections to earlier results and open problems}

\subsubsection*{Flat coalescence and exit clusters}
Basu, Busani and Ferrari proved the sharp cubic decay of the Airy$_1$ covariance \cite[Theorem~1.1]{BBF23}. Corollary~\ref{cor:app-flat} transfers the exponent and the constant to the flat coalescence probability, giving $\log p_{0,T}(r)\sim-r^3/(3T^2)$ as $rT^{-2/3}\to\infty$. In their exponential lattice model the two starts are separated by $u$ in the transverse coordinate of \cite[Theorem~2.7, (2.43)]{BBF23}. They bound the meeting probability by $\exp(-u^3/3+Cu^2)$ for $1<u\ll N^{1/14}$ and expect the leading term to be sharp. A matching lattice lower bound in this growing-$u$ regime requires estimates uniform in $N$ and $u$.

For a finite ordered family of starts, \eqref{eq:app-all} and \eqref{eq:app-count} give the probability of a single exit and the expected number of distinct exits. In the flat case \eqref{eq:app-cluster} gives the full joint survival function of the distances to the two ends of the cluster containing the origin. At small separations, combining \eqref{eq:app-airy} and \eqref{eq:app-short} at $T=1$ gives
\[
 \lim_{u\downarrow0}\frac{1+g_1'(u)}{u}=2^{2/3}\lambda_0 .
\]
Thus the one-sided curvature of the Airy$_1$ covariance measures the intensity of changes of the optimal exit. L\'opez and Pimentel prove a lower bound of order $n^{-2/3}$ for the survival probability of the tree of a prescribed flat lattice root and conjecture its asymptotic constant \cite[Theorem~2 and (2.24)]{LP17}. Our cluster is sampled at a fixed starting position. Relating its separation intensity to the prescribed-root survival constant also requires convergence of the discrete root counts under scaling \cite[Remark~2]{LP17}.

\subsubsection*{Time and position of first coalescence}
Rahman and Vir\'ag ask for the law of the space-time point of first coalescence \cite[Section~3.10]{RV25}. Pimentel's exponential-LPP duality expresses a coalescence-time probability as the probability of finding no exits in an interval \cite[Theorem~1]{Pim16}. For two zero-direction landscape geodesics from $(0,0)$ and $(r,0)$, \eqref{eq:app-time} gives the time marginal through the mixed derivative $\partial_q\partial_r V_{q,T}(r)$ at $q=1$. Boundary differentiation already gives a stationary meeting probability for finite-layer polymers at positive temperature \cite[Corollary~1.11]{RSZZ26}. For the landscape geodesics, Ganguly, Schweinsberg and Shuai give $\P(\tau_r>T)\le C rT^{-2/3}$ \cite[Proposition~3.11, (19)]{GSS26}. Here $\P(\tau_r>T)\sim\lambda_1 rT^{-2/3}$ at fixed $r$ as $T\to\infty$, with $\lambda_1$ the expected count in Proposition~\ref{prop:app-short}. The joint law with the position of first coalescence remains undetermined by \eqref{eq:app-time}.

The exit on the terminal line is a different spatial observable. For $q>0$, Pimentel's variance identity \cite[Theorem~2, (1.7)]{Pim22} and the absolute continuity in Lemma~\ref{lem:heights_exit} give the density of $-Y_{q,T}(0)$ as $v''(x)/(4q)$ for almost every $x$, where $v(x)=\Var(H_{q,T}(x))$. By \eqref{eq:app-scaling}, $\P(|Y_{q,T}(0)|\le R)\to0$ as $T\to\infty$ for each fixed $R<\infty$.

\subsubsection*{Exceptional endpoints and free exits}
Fix $a<b$. Let $D_{a,b}$ consist of the terminal positions reachable from $a$ and $b$ by a pair of geodesics disjoint at every time in $(0,T)$. Bates, Ganguly and Hammond prove $\dim_{\rm H}D_{a,b}=1/2$ almost surely \cite[Theorem~1.9(b)]{BGH22}. Writing $G_{a,b}(y)$ for the weight lost by imposing disjointness from $(a,0),(b,0)$ to the common endpoint $(y,T)$, we have
\[
 D_{a,b}\subseteq\{y:G_{a,b}(y)=0\},\qquad
 G_{a,b}(Y(a))>0\quad\text{a.s.}
\]
The inclusion follows by taking the disjoint geodesics as competitors. The strict inequality is Lemma~\ref{lem:strict_loss}. Hence $\P(Y(a)\in D_{a,b})=0$ even without conditioning on a common exit. This answers the selection question despite the dependence between the exit and the exceptional set. Dauvergne and Pankratz \cite[Theorem~1.1]{DP26} use the zero set and coordinate-section minima of a gap field for two paths with both endpoints in common to classify fixed-time geodesic networks. The gap above tests two distinct starts at an exit selected by the boundary reward.

\section{The finite-layer identity and its scaling}\label{sec:finite}

\subsection{The Brownian boundary}
Let $B_0,\ldots,B_m$ be independent two-sided standard Brownian motions pinned at zero. A path starting at $X$ takes one increment from each of the first $m$ rows. On reaching the last row at $P$, it pays its horizontal travel cost $P-X$ and receives the boundary reward $\alpha B_m(P)$. The two passage values are
\begin{align}
 G_m(X,P)&=\max_{X=u_0\le\cdots\le u_m=P}
                \sum_{j=0}^{m-1}\bigl[B_j(u_{j+1})-B_j(u_j)\bigr],\qquad P\ge X, \notag\\
 L(X)&=\max_{P\ge X}\{G_m(X,P)-(P-X)+\alpha B_m(P)\}.       \label{eq:finite-model}
\end{align}
For a fixed start the maximizing path is almost surely unique \cite[Lemma~2.4]{RSZZ26}. Write its jump positions in absolute coordinates as $t_j(X)=X+\Gamma_X(j)$, with $t_0(X)=X$ and exit $Q_m(X)=t_m(X)$ \cite[proof of Lemma~2.6, following (2.55)]{RSZZ26}. The two paths meet if their closed intervals intersect on some bulk row, that is, on the event
\[
 \mathcal M_m(Z,X)=\bigcup_{j=0}^{m-1}
 \{[t_j(Z),t_{j+1}(Z)]\cap[t_j(X),t_{j+1}(X)]\ne\varnothing\}.
\]
The finite-layer identity of \cite[Section~1.1 and Corollary~1.12]{RSZZ26} uses two-sided Brownian motions pinned at zero and allows arbitrary distinct starts $X,Z\in\R$. Its nonnegative coordinates are the displacements $\Gamma_X(j)=t_j(X)-X$. To convert the passage-time product and exit probability in that identity to our notation, observe that the increments of each row viewed from $X$ have the same law as those viewed from $0$. Hence
\[
 (L(X)-\alpha B_m(X),\,Q_m(X)-X)\ \overset d=\ (L(0),Q_m(0)),
 \qquad \E L(X)=\E L(0).
\]
In particular, $\P(Q_m(X)>0)=\P(Q_m(0)>-X)$, and the constant mean gives $\partial_X\E[L(Z)L(X)]=\partial_X\Cov(L(Z),L(X))$. With these substitutions the integration-by-parts identity \cite[Corollary~1.12]{RSZZ26} reads, for $Z\ne X$,
\begin{equation}\label{eq:finite_deriv}
 \partial_X\Cov(L(Z),L(X))
 =\tfrac12(1-\alpha^2)\sgn(Z-X)\P(\mathcal M_m(Z,X))
 +\alpha^2\bigl[\P(Q_m(X)>0)-\ind_{\{Z<X\}}\bigr].
\end{equation}

\subsection{The characteristic direction}
Take $k=\lfloor NT\rfloor$, $m=k+1$ and $T_N=k/N$. We center the exit window at the Brownian pin. A path from $-k$ has characteristic exit at $0$, so the start and exit coordinates are
\[
 X_N(a)=-k+2N^{2/3}a,\qquad P_N(y)=2N^{2/3}y,\qquad \ell=P_N(y)-X_N(a).
\]
The factor $2$ in the spatial coordinates fixes the coefficient of the limiting parabola. Indeed, for fixed $a,y$ and large $N$,
\[
 N^{-1/3}(2\sqrt{k\ell}-\ell-k)
 =-\frac{(\ell-k)^2}{N^{1/3}(\sqrt\ell+\sqrt k)^2}
 =-\frac{4N(y-a)^2}{(\sqrt\ell+\sqrt k)^2}\longrightarrow-\frac{(y-a)^2}{T}.
\]
We subtract the travel cost $\ell$ and the characteristic height $k$ from the bulk weight. Brownian scaling gives the boundary field in the same coordinates,
\begin{equation}\label{eq:center}
 S_N(a,y)=\frac{G_m(X_N(a),P_N(y))-\ell-k}{N^{1/3}},\qquad
 f_N(y)=\frac{\alpha B_m(2N^{2/3}y)}{N^{1/3}},\qquad f_N\overset d=f .
\end{equation}
Set $S_N=-\infty$ for $\ell<0$. Maximization in $y$ now yields
\[
 H_N(a)=\max_y(S_N(a,y)+f_N(y))=\frac{L(X_N(a))-k}{N^{1/3}},\qquad
 Y_N(a)=\frac{Q_m(X_N(a))}{2N^{2/3}}.
\]
Write $M_N(a,b)=\mathcal M_m(X_N(a),X_N(b))$. The two height factors contribute $N^{-2/3}$ to the covariance, and differentiating the starting coordinate contributes $X_N'(x)=2N^{2/3}$. Thus \eqref{eq:finite_deriv} gives, for $z\ne x$,
\begin{equation}\label{eq:scaled_deriv}
 \partial_x\Cov(H_N(z),H_N(x))
 =(1-\alpha^2)\sgn(z-x)\P(M_N(z,x))
 +2\alpha^2\bigl[\P(Y_N(x)>0)-\ind_{\{z<x\}}\bigr].
\end{equation}

\subsection{Fixed endpoints and a free exit}
To compare a pair of optimizers with a disjoint pair we need two fields. The one-path field is $S(a,u)=\mathcal L(a,0;u,T)$. For $a\le b$ and $u\le v$, write $S^{(2)}((a,b),(u,v))$ for the largest total length of paths from $a$ to $u$ and $b$ to $v$ that are disjoint on $(0,T)$. Its finite-layer counterpart $S_N^{(2)}$ uses essentially disjoint paths and the centering in \eqref{eq:center} for each path \cite[(6), Definition~1.1 and (7)]{DZ24}.

Subtracting the disjoint-pair weight from $S(a,u)+S(b,v)$ gives a continuous function of the ordered endpoints. At a common maximizing exit this difference is almost surely positive (Lemma~\ref{lem:strict_loss}). Its sign thus detects an interior meeting in the limit even when disjoint paths can approach a common segment.

For the coordinate convention of \cite[(11)]{DZ24}, set
$\widehat B_{-j}(s)=N^{-1/2}[B_j(Ns-k)-B_j(-k)]$. The finite endpoints in the new rows are
\[
 \left(\frac{X_N(a)+k}{N},0\right)=(2aN^{-1/3},0),\qquad
 \left(\frac{P_N(y)+k}{N},-k\right)=(T_N+2yN^{-1/3},-k).
\]
Dauvergne and Zhang prove joint convergence in law of the one- and two-path fields, uniformly on compact sets in the spatial endpoints and both time coordinates $s<t$ \cite[Theorem~1.5]{DZ24}. The spatial endpoints are weakly ordered, so these sets may include $u=v$. The limit is continuous and agrees with the disjoint variational field \cite[Theorems~1.4 and~1.6]{DZ24}. We may therefore evaluate at $(s,t)=(0,T_N)\to(0,T)$ and adjoin the independent boundary $f_N\overset d=f$ to obtain
\begin{equation}\label{eq:joint_conv}
 (S_N,S_N^{(2)},f_N)\ \Longrightarrow\ (S,S^{(2)},f)
 \qquad\text{locally uniformly in the endpoints.}
\end{equation}

\begin{lemma}[Exit tightness and joint convergence]\label{lem:heights_exit}
For each deterministic $a$, the maximum in \eqref{eq:height} has an almost surely unique location $Y(a)$, whose law is absolutely continuous with respect to Lebesgue measure. For any finite collection of deterministic starts $a_{i,N}\to a_i$, the convergence
\[
 \bigl(H_N(a_{i,N}),Y_N(a_{i,N})\bigr)_{i=1}^p
 \ \Longrightarrow\ \bigl(H(a_i),Y(a_i)\bigr)_{i=1}^p
\]
holds jointly with \eqref{eq:joint_conv}. For $1\le i,j\le p$,
\[
 \Cov\bigl(H_N(a_{i,N}),H_N(a_{j,N})\bigr)\longrightarrow
 \Cov\bigl(H(a_i),H(a_j)\bigr).
\]
\end{lemma}

\begin{proof}
\textbf{Step 1. Continuum localization and absolute continuity of the exit law.}
Fix a compact interval $I$ of starts and write $F_a(y)=S(a,y)+f(y)$. The parabolic bound of \cite[Corollary~10.7]{DOV22}, restricted to times $0,T$, gives
\[
 \sup_{a\in I}\left|S(a,y)+\frac{(y-a)^2}{T}\right|
 \le C_{I,T}\bigl(1+\log^2(2+|y|)\bigr),\qquad y\in\R,
\]
with $C_{I,T}<\infty$ almost surely. Since $f(y)=o(|y|)$ as $|y|\to\infty$, almost surely
\[
 \lim_{R\to\infty}\sup_{a\in I,\,|y|>R}F_a(y)=-\infty,
 \qquad \inf_{a\in I}F_a(0)>-\infty.
\]
All maxima are therefore attained in one random compact interval. On $J=[c,d]$, write $\delta_J g=(g(c+t)-g(c))_{0\le t\le d-c}$. For fixed $a$, the local Brownian comparison \cite[Theorem~1.2(2), (12), Remark~1.3(3)]{Dau23} and independence of $f$ give
\[
 \operatorname{Law}(\delta_J F_a)
 \ll \operatorname{Law}\bigl((\sqrt{2(1+\alpha^2)}\,W(t))_{0\le t\le d-c}\bigr),
\]
where $W$ is standard Brownian motion. Its maximum on $[0,d-c]$ has a unique location with the arcsine law. Thus on each $[-n,n]$, $n\ge1$, $F_a$ has an almost surely unique maximizer $Y^{(n)}$ with an absolutely continuous law. Localization gives uniqueness of $Y(a)$, and for every Lebesgue-null set $A$,
\[
 \P(Y(a)\in A)\le\sum_{n\ge1}\P(Y^{(n)}\in A)=0 .
\]

\noindent\textbf{Step 2. Linear tilt and random-matrix edge bounds.}
We add a linear reward $hy$ to the bulk weight. Let $G_m^B$ denote passage in rows $B=(B_0,\ldots,B_{m-1})$ and set
\[
 S_N^h(a)=\sup_y\{S_N(a,y)+hy\},\qquad
 A_m(B)=\sup_{t\ge0}\{G_m^B(0,t)-t\},\qquad
 Z_N=\frac{A_m(B)-k}{N^{1/3}}.
\]
For fixed $h$ and large $N$, $\theta=1-h/(2N^{1/3})>0$. In the rows $\widetilde B_j(s)=\theta B_j(s/\theta^2)$, rescaling every jump coordinate gives
\[
 G_m^{\widetilde B}(0,s)=\theta G_m^B(0,s/\theta^2),\qquad
 A_m(\widetilde B)=\theta\sup_{t\ge0}\{G_m^B(0,t)-\theta t\}.
\]
These equalities hold pathwise, and $\widetilde B\overset d=B$. Set $y=a+(t-k)/(2N^{2/3})$, so that $P_N(y)=X_N(a)+t$. For $t\ge0$, \eqref{eq:center} gives
\[
 S_N(a,y)+hy
 =ha-\frac{hk}{2N^{2/3}}-\frac{k}{N^{1/3}}
   +\frac{G_m(X_N(a),X_N(a)+t)-\theta t}{N^{1/3}}.
\]
Taking the supremum over $t\ge0$ and translating the row increments to the origin gives
\begin{equation}\label{eq:tilt_bound}
 S_N^h(a)\overset d=ha-\frac{hk}{2N^{2/3}}+
 \frac{A_m(B)/\theta-k}{N^{1/3}}
 =ha+\frac{kh^2}{4N\theta}+\frac{Z_N}{\theta}.
\end{equation}
In the last equality, the deterministic correction is
$kN^{-1/3}(\theta^{-1}-1-(1-\theta))=kN^{-1/3}(1-\theta)^2/\theta$.

We need a fourth-moment bound for $Z_N$ and for the weight of the competitor $y=a$. Let $\mathcal H_m$ be standard Hermitian Brownian motion and $\lambda_H=\lambda_{\max}(\mathcal H_m(1))$. The Brownian last-passage representation \cite[(9) and Theorem~7]{OY02}, with signs reversed, and Brownian scaling give
\begin{equation}\label{eq:hermitian}
 (G_m^B(0,t))_{t\ge0}\overset d=(\lambda_{\max}(\mathcal H_m(t)))_{t\ge0},
 \qquad G_m(X,X+k)\overset d=\sqrt{k}\,\lambda_H .
\end{equation}
For an $(m+1)\times m$ matrix $D$ of independent standard real Gaussians, put $\lambda_W=\lambda_{\max}(D^{\mathsf T}D)$. The identity of \cite[Section~2, (5)]{FW20} gives
\begin{equation}\label{eq:wishart}
 4A_m(B)\overset d=4\sup_{t\ge0}\{\lambda_{\max}(\mathcal H_m(t))-t\}\overset d=\lambda_W .
\end{equation}
The Wishart edge has center $\mu_m=(\sqrt{m+1}+\sqrt m)^2$ and scale $\sigma_m=(\sqrt{m+1}+\sqrt m)^{4/3}/[m(m+1)]^{1/6}$. In the moment estimates of \cite[(1.3), (1.5) and the paragraph following Corollary~3]{LR10}, take $p=4$, $\beta=2$ for Hermite, and $(\beta,\kappa)=(1,m+1)$ for Laguerre. Then
\begin{equation}\label{eq:edge}
 \sup_{m\ge1}\left\{
  \E|m^{1/6}(\lambda_H-2\sqrt m)|^4+
  \E|(\lambda_W-\mu_m)/\sigma_m|^4\right\}<\infty .
\end{equation}
Write $\xi_H=m^{1/6}(\lambda_H-2\sqrt m)$ and $\xi_W=(\lambda_W-\mu_m)/\sigma_m$. In our coordinates \eqref{eq:hermitian} and \eqref{eq:wishart} read
\[
 \begin{aligned}
 S_N(a,a)&\overset d=\frac{\sqrt{k}\,m^{-1/6}}{N^{1/3}}\xi_H
                         +\frac{2\sqrt{km}-2k}{N^{1/3}},\\
 Z_N&\overset d=\frac{\sigma_m}{4N^{1/3}}\xi_W+\frac{\mu_m-4k}{4N^{1/3}}.
 \end{aligned}
\]
Here $m=k+1$, $k/N\to T$ and $\sigma_m\asymp m^{1/3}$. The deterministic centerings and the resulting bound from \eqref{eq:edge} are
\[
 2\sqrt{km}-2k=O(1),\qquad \mu_m-4k=O(1),\qquad
 \sup_{N\ge N_0,\,a\in I}\bigl(\|S_N(a,a)\|_4+\|Z_N\|_4\bigr)<\infty .
\]

\noindent\textbf{Step 3. Uniform fourth moments of heights and exits.}
Set $K_N=\sup_y\{f_N(y)-|y|\}$ and $K=\sup_y\{f(y)-|y|\}$. By \eqref{eq:center}, $K_N\overset d=K$ for every $N$. For $\alpha\ne0$, the exponential martingale $\exp((f(t)-t)/\alpha^2)$ on $t\ge0$ and its counterpart on the negative half-line give $\P(K>u)\le 2e^{-u/\alpha^2}$ for $u\ge0$. For $\alpha=0$, $K=0$. In particular $\sup_N\|K_N\|_4<\infty$. Write $H_N^h(a)=\sup_y\{S_N(a,y)+f_N(y)+hy\}$. Evaluating at $y=a$ and using $f_N(y)\le K_N+|y|$ gives
\[
 S_N(a,a)+f_N(a)\le H_N(a)=H_N^0(a),\qquad
 H_N^h(a)\le K_N+\max_{\epsilon=\pm1}S_N^{h+\epsilon}(a).
\]
At $y=Y_N(a)$, one of the slopes $h=\pm1$ adds $|Y_N(a)|$, so
\[
 |Y_N(a)|\le\max_{\epsilon=\pm1}H_N^\epsilon(a)-H_N^0(a).
\]
Since $\sup_{N,\,a\in I}\|f_N(a)\|_4=\sup_{a\in I}\|f(a)\|_4<\infty$, the preceding bounds and \eqref{eq:tilt_bound} at $h=0,\pm1,\pm2$ give
\begin{equation}\label{eq:heights_exit_moment}
 \sup_{N\ge N_0,\,a\in I}
 \E\bigl[|H_N(a)|^4+|Y_N(a)|^4\bigr]\le C_I<\infty .
\end{equation}

\noindent\textbf{Step 4. Joint convergence and convergence of covariances.}
Choose $I$ to contain all the starts. By \eqref{eq:heights_exit_moment},
$\P(\max_i|Y_N(a_{i,N})|>R)\le pC_I R^{-4}$ for $N\ge N_0$.
Together with \eqref{eq:joint_conv}, this bound gives a further subsequence of every subsequence on which the fields and exits converge jointly in law. For the fields we use a countable exhaustion by compact endpoint sets, each feasible for all sufficiently large $N$. Relabel this further subsequence by $N$. The Skorokhod representation \cite[Theorem~6.7]{Bil99} gives a coupling on which the fields converge locally uniformly and $Y_N(a_{i,N})\to y_i$ for every $i$, almost surely. Write $F_{i,N}(y)=S_N(a_{i,N},y)+f_N(y)$. Optimality gives, simultaneously for $q\in\mathbb Q$,
\[
 F_{a_i}(y_i)=\lim_N F_{i,N}(Y_N(a_{i,N}))
             \ge\lim_N F_{i,N}(q)=F_{a_i}(q).
\]
Continuity and uniqueness force $y_i=Y(a_i)$ and $H_N(a_{i,N})\to H(a_i)$ in this coupling. This identifies the joint limit along every subsequence. Finally, with $H_{i,N}=H_N(a_{i,N})$ and $H_i=H(a_i)$, \eqref{eq:heights_exit_moment} implies
\[
 \sup_{N\ge N_0}\|H_{i,N}H_{j,N}\|_2
 \le \sup_{N\ge N_0}\|H_{i,N}\|_4\|H_{j,N}\|_4<\infty .
\]
Uniform integrability gives $\E H_{i,N}\to \E H_i$ and $\E[H_{i,N}H_{j,N}]\to \E[H_iH_j]$, hence the asserted convergence of covariances.
\end{proof}

\section{Coalescence at the maximizing exit}\label{sec:coalescence}

Suppose two geodesics meet before the terminal line. At their meeting point either geodesic can take the other's continuation. We first establish the uniqueness needed for this exchange, then consider the converse question at a common exit.

\begin{lemma}[Free-endpoint geodesics]\label{lem:free_geodesics}
For each deterministic $a$, the maximizing path $\pi_a$ in \eqref{eq:height} is almost surely unique. For fixed $a<b$, almost surely $Y(a)\le Y(b)$ and
\begin{equation}\label{eq:tails}
 \pi_a(r)=\pi_b(r)\quad\Longrightarrow\quad
 \pi_a|_{[r,T]}=\pi_b|_{[r,T]},\qquad 0<r<T .
\end{equation}
\end{lemma}

\begin{proof}
Lemma~\ref{lem:heights_exit} gives a unique exit $Y(a)$. Since the landscape admits geodesics simultaneously for all endpoint pairs \cite[Theorem~1.7, first assertion]{DZ24}, there is a geodesic from $(a,0)$ to $(Y(a),T)$.

Fix a deterministic time $0<r<T$ and set
\[
 A_r(v)=\mathcal L(a,0;v,r),\qquad
 F_r(v)=\max_y\{\mathcal L(v,r;y,T)+f(y)\}.
\]
The parabolic bound \cite[Corollary~10.7]{DOV22} confines the exits for $v$ in a compact interval to a common random compact set. Thus $F_r$ is continuous. The independent-increment and metric composition properties \cite[Definition~10.1(II) and Lemma~10.8, (10.12)]{DOV22} give, for every maximizing path $\pi$,
\[
 A_r\ \text{independent of }F_r,\qquad
 H(a)=\max_v(A_r+F_r)(v)=(A_r+F_r)(\pi(r)).
\]

We show that the position $\pi(r)$ is unique. Condition on $F_r=g$. On a compact interval the recentered $A_r$ is absolutely continuous in law with respect to Brownian motion of variance rate $2$ \cite[Theorem~1.2(2)]{Dau23}, \cite[Lemma~10.2(5)]{DOV22}. Let $I=[p,q]$ and $J=[u,v]$ have rational endpoints with $q<u$, and write $M_I(h)=\max_I h$. For a comparison Brownian motion $W$ starting at $p$,
\[
 \begin{aligned}
 M_J(W+g)-M_I(W+g)&=W(u)-W(q)+K,\\
 K&=W(q)-M_I(W+g)+M_J(W-W(u)+g).
 \end{aligned}
\]
The first two terms in $K$ depend only on $W|_{[p,q]}$, and the last depends only on the increments after $u$. Hence $W(u)-W(q)$ is independent of $K$ and has law $N(0,2(u-q))$. The difference of the two maxima therefore has no atom at zero. Absolute continuity gives
\[
 \P\{M_I(W+g)=M_J(W+g)\}=0
 \quad\Longrightarrow\quad \P\{M_I(A_r+g)=M_J(A_r+g)\}=0 .
\]
Any two distinct maximizers would give a tie on some such pair of intervals. Almost surely all maximizing paths therefore agree at every rational $r\in(0,T)$, and continuity gives uniqueness on $[0,T]$.

At a meeting $(v,r)$, simultaneous metric composition gives
\[
 \mathcal L(v,r;Y(a),T)+f(Y(a))
 =F_r(v)=\mathcal L(v,r;Y(b),T)+f(Y(b)).
\]
Exchanging the two continuations leaves both paths optimal, so uniqueness gives \eqref{eq:tails}. If $Y(a)>Y(b)$, the continuous difference $\pi_a-\pi_b$ is negative at time $0$ and positive at time $T$. It must vanish at some $r\in(0,T)$, where \eqref{eq:tails} forces $Y(a)=Y(b)$.
\end{proof}

\begin{lemma}[Strict loss at the maximizing exit]\label{lem:strict_loss}
For fixed $a<b$, let
\[
 G_{a,b}(y)=S(a,y)+S(b,y)-S^{(2)}((a,b),(y,y)).
\]
Then $G_{a,b}(Y(a))>0$ almost surely, and
$M(a,b)=\{Y(a)=Y(b)\}$ almost surely.
\end{lemma}

\begin{proof}
By landscape scaling \cite[Lemma~10.2(5)]{DOV22} and Brownian scaling we may take $T=1$. Dauvergne's two-line representation \cite[(12)]{Dau23} gives
\[
 W_1(y)=S(a,y),\qquad W_1(y)+W_2(y)=S^{(2)}((a,b),(y,y)).
\]
Set $D=W_2-W_1$. The first equality below is the running-maximum formula \cite[Theorem~1.4(1)]{Dau23}. Substitution in the definition of $G_{a,b}$ gives the second.
\begin{equation}\label{eq:record}
 S(b,y)-S(a,y)=\sup_{w\le y}D(w),\qquad
 G_{a,b}(y)=\sup_{w\le y}D(w)-D(y).
\end{equation}
The exit $Y(a)$ maximizes $U=W_1+f$. By \eqref{eq:record} it suffices to find a point to its left where $D$ takes a larger value.

Fix $I=[p,q]$ with rational endpoints $p<q$, and write $g_I(y)=g(y)-g(p)$ for $y\in I$. By local Brownian absolute continuity \cite[Theorem~1.2(2)]{Dau23} and independence of $f$, the law of $((W_1)_I,(W_2)_I,f_I)$ is absolutely continuous with respect to the law of three independent Brownian motions on $I$, started at zero at $p$, with variance rates $2,2,2\alpha^2$ (the third path is zero when $\alpha=0$). Under this comparison law, write $\widehat U$ and $\widehat D$ for the sum of the first and third paths and the difference of the second and first paths, respectively. Their covariance matrix per unit time is
\[
 \begin{pmatrix}2(1+\alpha^2)&-2\\-2&4\end{pmatrix}.
\]
Set $V=\widehat D+\widehat U/(1+\alpha^2)$. For $s,t\in I$,
\[
 \Cov(\widehat U(s),V(t))
 =(\min\{s,t\}-p)\left(-2+\frac{2(1+\alpha^2)}{1+\alpha^2}\right)=0 .
\]
Joint Gaussianity gives independence of $V$ and $\widehat U$. The variance rate of $V$ is $4-2/(1+\alpha^2)>0$.

Let $\widehat y_I$ maximize $\widehat U$ on $I$. For each fixed $y\in I^\circ$, the backward increments $V(y-h)-V(y)$ have Brownian law in $h\in[0,y-p]$. The law of the iterated logarithm at zero gives
\[
 \limsup_{h\downarrow0}\frac{V(y-h)-V(y)}{\sqrt h}=+\infty\qquad\text{almost surely}.
\]
Conditional on $\widehat U$, its unique maximizer $\widehat y_I$ is fixed and $V$ retains its Brownian law. On $\{\widehat y_I\in I^\circ\}$ there are therefore points $w\uparrow\widehat y_I$ for which
\[
 \widehat D(w)-\widehat D(\widehat y_I)
 =\frac{\widehat U(\widehat y_I)-\widehat U(w)}{1+\alpha^2}
   +V(w)-V(\widehat y_I)>0 .
\]

The maximizer on $I$ and the differences of $D$ on $I$ depend only on the recentered paths. Absolute continuity transfers the preceding inequality to $D$ at the maximizer of $U$ on $I$ whenever it lies in $I^\circ$. On $\{Y(a)\in I^\circ\}$ that maximizer is $Y(a)$, so \eqref{eq:record} gives
\[
 \P\{Y(a)\in I^\circ,\ G_{a,b}(Y(a))=0\}=0 .
\]
Taking the union over all rational $p<q$ gives $\P\{G_{a,b}(Y(a))=0\}=0$, since the interiors of these intervals cover $\R$. Now $G_{a,b}\ge0$, which proves the strict inequality.

If two geodesics to a common exit $y$ were disjoint on $(0,T)$, their combined length $S(a,y)+S(b,y)$ would be admissible in the disjoint variational formula \cite[Definition~1.1, (7)]{DZ24}. This would give $G_{a,b}(y)=0$. Hence a common exit forces an interior meeting. The reverse implication is \eqref{eq:tails}.
\end{proof}

\begin{proposition}[Convergence of meeting probabilities]\label{prop:meeting_conv}
If $(a_N,b_N)\to(a,b)$ deterministically with $a<b$, then
$\P(M_N(a_N,b_N))\to\P(M(a,b))$.
\end{proposition}

\begin{proof}
Write $u_N=Y_N(a_N)$, $v_N=Y_N(b_N)$ and abbreviate the meeting event by $M_N$. For all large $N$, $a_N<b_N$. On $M_N$ the two paths have a common point on a bulk row. Exchanging tails at that point preserves the order of the jump positions on each path and the sum of the two weights, including the terminal rewards. Each exchanged path must therefore attain its starting optimum. Finite-path uniqueness \cite[Lemma~2.4]{RSZZ26} forces the tails to agree, so $M_N\subseteq\{u_N=v_N\}$. Since the starts are ordered, reversed exits would force an intersection on a bulk row and hence equal exits. Thus $u_N\le v_N$. Their separation and the cost of disjointness at these two exits are
\[
 \begin{aligned}
 d_N&=v_N-u_N,\\
 \Delta_N&=S_N(a_N,u_N)+S_N(b_N,v_N)
       -S_N^{(2)}((a_N,b_N),(u_N,v_N)).
 \end{aligned}
\]
Write $X_a=X_N(a_N)$ and $X_b=X_N(b_N)$. On $M_N^c$ the closed intervals on each bulk row are disjoint. Starting from $t_0(X_a)<t_0(X_b)$ and proceeding row by row gives
\[
 t_{j+1}(X_a)<t_j(X_b)\le t_{j+1}(X_b),\qquad 0\le j<m .
\]
Thus the two paths are ordered and occupy different rows at every horizontal coordinate interior to both domains. They are essentially disjoint in the sense of \cite[Section~2.1]{DZ24} and attain the sum of the two individual optima, so $\Delta_N=0$. It follows that
\begin{equation}\label{eq:inclusion}
 \{\Delta_N>0\}\subseteq M_N\subseteq\{d_N=0\},\qquad
 d_N,\Delta_N\ge0,\quad d_N\Delta_N=0 .
\end{equation}
Define $d,\Delta$ by the same formulas at $u=Y(a)$, $v=Y(b)$. The extended landscape is continuous on the ordered endpoint space, including $u=v$, and agrees there with the disjoint variational field \cite[Theorems~1.4 and~1.6]{DZ24}. Evaluation at the exits is therefore continuous for the joint convergence in \eqref{eq:joint_conv} and Lemma~\ref{lem:heights_exit}. The continuous mapping theorem gives $(d_N,\Delta_N)\Rightarrow(d,\Delta)$.

By \eqref{eq:inclusion}, these laws are supported on $E=\{(s,t)\in[0,\infty)^2:st=0\}$. The lower and upper events in \eqref{eq:inclusion} correspond to the open and closed vertical rays in $E$,
\[
 O=E\cap\{t>0\},\qquad F=E\cap\{s=0\}=\overline O^{\,E},\qquad
 \partial_E O=F\setminus O=\{(0,0)\}.
\]
At $d=0$ the exits coincide and $\Delta=G_{a,b}(u)>0$ by Lemma~\ref{lem:strict_loss}. Consequently the limit law $\mu$ satisfies
\begin{equation}\label{eq:null-boundary}
 \mu(\partial_E O)=\P(d=0,\Delta=0)=0,\qquad
 \mu(O)=\mu(F)=\P(M(a,b)).
\end{equation}
Let $\mu_N$ be the law of $(d_N,\Delta_N)$. Applying the Portmanteau bounds \cite[Theorem~2.1(iii),(iv)]{Bil99} to $O$ and $F$, with \eqref{eq:inclusion} between them, gives
\begin{equation}\label{eq:squeeze}
 \begin{aligned}
 \mu(O)&\le\liminf_N\mu_N(O)\le\liminf_N\P(M_N)\\
       &\le\limsup_N\P(M_N)\le\limsup_N\mu_N(F)\le\mu(F).
 \end{aligned}
\end{equation}
Both endpoints equal $\P(M(a,b))$ by \eqref{eq:null-boundary}.
\end{proof}

\begin{proof}[Proof of Theorem~\ref{thm:main} and Corollary~\ref{cor:integral}]
Write $R_N(z,x)$ for the right side of \eqref{eq:scaled_deriv} and $R(z,x)$ for the same expression with the limiting exits and meeting probabilities. For every deterministic sequence $(z_N,x_N)\to(z,x)$, Lemma~\ref{lem:heights_exit} gives
\[
 \Cov(H_N(z_N),H_N(x_N))\longrightarrow\Cov(H(z),H(x)).
\]
For $z\ne x$, Proposition~\ref{prop:meeting_conv} gives convergence of the meeting probability. Since $\P(Y(x)=0)=0$, the exit probability also converges. Hence $R_N(z_N,x_N)\to R(z,x)$ off the diagonal.

Consequently $(z,x)\mapsto\Cov(H(z),H(x))$ is continuous on $\R^2$, and $R(z,x)$ is continuous on $\{(z,x):z\ne x\}$. The corresponding convergences are locally uniform on these domains. Indeed, suppose $g_N(w_N)\to g(w)$ whenever $w_N\to w$. Given $w_j\to w$, choose $N_j\uparrow\infty$ with $|g_{N_j}(w_j)-g(w_j)|<1/j$. Then
\[
 |g(w_j)-g(w)|\le j^{-1}+|g_{N_j}(w_j)-g(w)|\longrightarrow0 .
\]
On a compact set, a subsequence violating uniform convergence would have points $w_{N_j}\to w$ along a further subsequence. There $g_{N_j}(w_{N_j})-g(w_{N_j})\to0$.

Fix $z\notin[u,v]$. By \eqref{eq:scaled_deriv}, the finite-layer covariance is differentiable in $x$ on $[u,v]$ with $|R_N(z,x)|\le |1-\alpha^2|+2\alpha^2$. The mean value theorem makes it Lipschitz, hence absolutely continuous. The fundamental theorem of calculus and dominated convergence give
\[
 \begin{aligned}
 \Cov\bigl(H(z),H(v)-H(u)\bigr)
 &=\lim_{N\to\infty}\Cov\bigl(H_N(z),H_N(v)-H_N(u)\bigr)\\
 &=\lim_{N\to\infty}\int_u^v R_N(z,x)\,dx
   =\int_u^v R(z,x)\,dx .
 \end{aligned}
\]
By continuity of $R$, the derivative with respect to the moving start is $R(z,x)$. Taking $(z,x)=(a,b)$ gives the first line of \eqref{eq:cov_deriv}. Taking $(z,x)=(b,a)$ and using symmetry of covariance gives the second. Finally, integrate the first line from $a+\varepsilon$ to $b$. For $x\in I=[a,a+1]$, the localization in the proof of Lemma~\ref{lem:heights_exit} restricts the maximum in \eqref{eq:height} to a common random compact interval. Uniform continuity of $(x,y)\mapsto\mathcal L(x,0;y,T)$ on the resulting compact rectangle gives $H(a+\varepsilon)\to H(a)$ almost surely. For each $x\in I$, Lemma~\ref{lem:heights_exit} gives $H_N(x)\Rightarrow H(x)$. The lower semicontinuity of fourth moments and \eqref{eq:heights_exit_moment} give
\[
 \sup_{x\in I}\E|H(x)|^4
 \le \sup_{x\in I}\liminf_{N\to\infty}\E|H_N(x)|^4\le C_I .
\]
Thus the squared increments are uniformly integrable and $H(a+\varepsilon)\to H(a)$ in $L^2$. By Cauchy--Schwarz,
\[
 \bigl|\Cov(H(a),H(a+\varepsilon))-\Var(H(a))\bigr|
 \le \sqrt{\Var(H(a))}\,\|H(a+\varepsilon)-H(a)\|_2\longrightarrow0 .
\]
Bounded convergence in the integral gives \eqref{eq:cov_integral}.
\end{proof}

\end{document}